\documentclass[a4,12pt]{article}

\usepackage{amsmath,amssymb}

\newtheorem{thm}{Theorem}
\newtheorem{prop}{Proposition}
\newtheorem{lem}{Lemma}

\newtheorem{rmk}{Remark}
\newtheorem{dfn}{Definition}
\newtheorem{expl}{Example}

\title{An inequality for the $\Delta$-genus  of toric varieties}
\author{Shoetsu Ogata and  Riki Tabei}

\begin{document}
\maketitle

\begin{abstract}
For a polarized toric variety $(X,L)$ of dimension $n\le4$, we give a lower bound of the $\Delta$-genus
by using the vanishing number of adjoint bundle of a multiple of $L$.
We show that 
for a polarized toric variety of dimension $n$ with nonvanishing adjoint bundle, the $\Delta$-genus
is more than or equal to $n-1$.
\end{abstract}

\section*{Introduction}

A pair $(X,L)$ of a projective variety $X$ and an ample line bundle $L$ on $X$ is called a {\it polarized variety}.
For a polarize variety $(X,L)$ of dimension $n$, Fujita (see \cite{Fj}) defined the $\Delta$-genus as
$$
\Delta(X,L) :=n + L^n -\dim \Gamma(X, L).
$$

When $X$ is a toric variety, we call $(X,L)$ a polarized toric variety.

\begin{thm}\label{t:A}
Let $(X,L)$ be a polarized toric variety of dimension $n\le4$.
Then we have
\begin{equation}\label{e:a1}
\Delta(X,L) \ge n-\min\{k\ge1:\ \Gamma(X,L^{\otimes k}\otimes \omega_X)\not=0\},
\end{equation}
where $\omega_X$ is the dualizing sheaf of $X$.
\end{thm}

The statement of Theorem~\ref{t:A} is from Theorem~\ref{th:BN1}, Lemma~\ref{lem:5.1}, Lemma~\ref{lem:5.2}
and Proposition~\ref{p:5}. In order to unify inequalities into this form, 
we have to restrict dimension $n\le4$ because of Example~\ref{ex:3.2}. 

\begin{thm}\label{t:C}
Let $(X,L)$ be a polarized toric variety of dimension $n\ge3$.
Assume that $\Gamma(X,L^{\otimes n-2}\otimes\omega_X)=0$ and $\Gamma(X,L^{\otimes n-1}\otimes\omega_X)\not=0$.
Then $\Delta(X,L)\ge1$.  Moreover the equality holds if and only if $(X,L)$ is a Gorenstein toric Del Pezzo variety, that is,
$X$ is a Gorenstein toric variety with $L^{n-1}\cong \omega_X^\vee$.
\end{thm}

This theorem is given as Lemma~\ref{lem:5.1}.

\begin{rmk}
The assumption in Theorem~\ref{t:C} is necessary. See Example~\ref{ex:3.2}.  The polytope defines a polarized Gorenstein toric Fano
variety $(X,L)$ not Del Pezzo with $\Delta(X,L)=1$.
\end{rmk}

\begin{thm}\label{t:B}
Let $(X,L)$ be a polarized toric variety of dimension $n$ with $\Gamma(X,L\otimes\omega_X)\not=0$.
Then we have
$$
\Delta(X,L) \ge n-1.
$$
Moreover, when $3\le n$, the equality holds if and only if $(X,L)$ is the Gorenstein toric Fano variety which is
a mirror to $(\mathbb{P}^n, {\cal O}(n+1))$.
\end{thm}

This theorem is given as Proposition~\ref{p:5}.

\bigskip

This research begins from the work of Tabei\cite{Tb}.  We give his main result in the last section.

\section{An origin of the inequality}

Let $X \subset \mathbb{P}^N$ be a  projective variety of dimension $n$ and ${\cal I}_X$ the ideal sheaf of $X$.
The {\it Castelnuovo-Mumford regularity} $\mbox{\rm reg}(X)$ of $X$ is defined by
$$
\mbox{\rm reg}(X):= \min\{r\ge1:\ H^i(\mathbb{P}^N, {\cal I}_X(r-i))=0 \quad \mbox{for all $i>0$}\}.
$$
When $X$ is nondegenerate and irreducible, Eisenbud and Goto conjectured 
\begin{equation}\label{e:1.1}
\mbox{\rm reg}(X) \le \mbox{\rm deg} X -\mbox{\rm codim} X +1.
\end{equation}

The conjecture has been proved true in several cases. In 2018,
Peeva and McCullough gave a counterexample to this inequality. 

Here we consider the situation that $X$ is a toric variety and an embedding given by global sections
of a very ample line bundle $L$, that is, $\Phi _L: X \hookrightarrow \mathbb{P}(\Gamma(X,L)^*)$.
In this case, $\mbox{\rm reg}(X)$ is determined by two values:
Set the bound of $k$-normality as
\begin{equation}\label{e:1.2}
\kappa(X):=\min\{k_0\ge1:\ \Gamma(X,L)^{\otimes k} \to \Gamma(X,L^{\otimes k})\quad \mbox{is surjective for all $k\ge k_0$}\}.
\end{equation}
And set
\begin{equation}\label{e:1.3}
\lambda(X,L):=\min\{k\ge1:\ \Gamma(X, L^{\otimes k}\otimes\omega_X)\not=0\}.
\end{equation}
Then we have
$$
{\rm reg}(X)=\max\{\kappa(X)+1, \dim X+2-\lambda(X,L)\}.
$$
From $\mbox{\rm reg}(X)\ge \dim X+2-\lambda(X,L)$, we have a weaker inequality
\begin{equation}\label{e:1.4}
\dim \Gamma(X,L) -\deg X \le \lambda(X,L).
\end{equation}
If  $\dim X=n$, then we write $\deg X=L^n$, and  all three values in the inequality (\ref{e:1.4}) can be defined for any ample
line bundle $L$.

Thus for a polarized toric variety $(X,L)$ of dimension $n$, we may consider an inequality
\begin{equation}\label{e:d1}
\Delta(X,L) \ge n-\lambda(X,L)
\end{equation}

\subsection{Polarized toric variety}

For a proof of theorems we use combinatrics of lattice polytopes. So we recall basic notions of toric varieties and lattice polytopes.
See, for example, Oda's book\cite{Od} or Fulton's book\cite{Ft}.
Let $M$ be a free abelian group of rank $n$.  Denote $M_{\mathbb{R}}\cong \mathbb{R}^n$ the extension of coefficients of $M$
to real numbers.
Let $\mathbb{C}[M]$ be the group algebra over complex numbers, and $T=\mbox{\rm Spec}\ \mathbb{C}[M]$ the algebraic torus
of dimension $n$. The character group $\mbox{\rm Hom}(T,\mathbb{C}^\times)$ of $T$ is isomorphic to $M$.  For a lattice point
$m\in M$, the corresponding character is written by $\chi^m$.
A toric variety $X$ of dimension $n$ is a normal algebraic variety with an algebraic action of the algebraic torus $T$ of dimension $n$
such that $X$ has an open orbit isomorphic to $T$ and if we identify the open orbit with $T$ the action of $T$ is compatible with
the multiplication of $T$.

A lattice polytope $P$ in $M_{\mathbb{R}}$ is a convex hull of finite elements of $M$ in $M_{\mathbb{R}}$.
A polarized toric variety $(X,L)$ of dimension $n$ corresponds to a lattice polytope $P\subset M_{\mathbb{R}}$ of dimension $n$.
The correspondence implies the equalities
\begin{eqnarray*}
\Gamma(X,L) &=& \bigoplus_{m\in P\cap M} \mathbb{C} \chi^m,\\
\Gamma(X,L\otimes\omega_X) &=& \bigoplus_{m\in \mbox{\scriptsize int}P\cap M} \mathbb{C} \chi^m.
\end{eqnarray*}
A twist $L^{\otimes k}$ corresponds to the multiple $kP$.  The condition that the multiplication map $\Gamma(X,L)^{\otimes k} \to
\Gamma(X,L^{\otimes k})$ is surjective is equivalent to the equality
\begin{equation}\label{e:1.6}
\overbrace{(P\cap M) + \dots +(P\cap M)}^k = (kP)\cap M.
\end{equation}
A lattice polytope $P\subset M_{\mathbb{R}}$ is {\it normal} if the equation (\ref{e:1.6}) holds for all $k\ge1$.
Here we define {\it empty depth} of a lattice polytope $P$ as
$$
e(P) : = \begin{cases} \max\{k\ge1:\ \mbox{int}(kP)\cap M=\emptyset\} & \mbox{if $\mbox{int}P\cap M =\emptyset$},\\
0& \mbox{if $\mbox{int}P\cap M \not=\emptyset$}. \end{cases}
$$
Then we see $\lambda(X,L)=e(P)-1$.
A lattice polytope $P$ of dimension $n$ has the normalized volume $n! \mbox{vol}(P)$, which we write as $v(P)$.
Then the degree $L^n$ coincides with $v(P)$.  The $\Delta$-genus $\Delta(X,L)$ is written as
$$
\Delta(P) = \dim P +v(P) -\sharp(P\cap M).
$$
The inequality (\ref{e:d1}) is equivalent to
\begin{equation}\label{e:p1}
\sharp(P\cap M) -v(P) \le e(P) +1.
\end{equation}

\subsection{The case of dimension two}

Let $P$ be a lattice polytope of dimension two. We know Pick's formula concerning  the area and the number of
lattice points in $P$ (see, for example, \cite[p.101]{Od}).

\begin{thm}[Pick]
Let $P\subset M_{\mathbb{R}}$ be a lattice polygon, that is , a lattice polytope of dimension two.
Then we have
\begin{equation}\label{e:2.1}
\mbox{\rm vol}(P) =\sharp (P\cap M) -\frac12\sharp(\partial P\cap M)-1.
\end{equation}
\end{thm}

Since $v(P)=2\mbox{vol}(P)$, we write the equality (\ref{e:2.1}) as
$$
\sharp(P\cap M) -v(P) = 2-\sharp(\mbox{int}P\cap M).
$$

\begin{thm}[Tabei\cite{Tb}]
Let $P\subset M_{\mathbb{R}}$ be a lattice polygon.  Then we have
$$
\sharp(P\cap M) -v(P) \le e(P) +1.
$$
Moreover, the equality holds if and only if $P$ is not basic and $\sharp(\mbox{\rm int}P\cap M)\le1$.
\end{thm}

\begin{rmk}
Since a lattice polygon is always normal, the inequality {\rm (\ref{e:d1})} is equivalent to {\rm (\ref{e:1.1})}.
\end{rmk}

\begin{prop}[Tabei\cite{Tb}, Koelman\cite{Ko}]
A lattice polygon $P$ with $e(P)=1$ is isomorphic to
$$
P=\mbox{\rm conv}\{0,(1,0),(0,a),(1,b)\}, \quad\mbox{where $a\ge1, b\ge0$ and $a+b\ge2$}
$$
or
$$
P=\mbox{\rm conv}\{0, (0,2),(2,0)\}.
$$
\end{prop}

\begin{rmk}
A lattice polygon $P$ with $\sharp(\mbox{int}P\cap M)=1$ is called "Fano polygon".
Fano polygons are classified and there are 16 polygons up to isomorphism.
See {\rm \cite{Ks} or \cite{KN}}.
\end{rmk}

\section{Empty lattice simplices}

A lattice polytope $P\subset M_{\mathbb{R}}$ of dimension $n$ is called {\it empty lattice simplex}
if $\sharp(P\cap M)=n+1$.
Here we consider $e(P)$ and $k$-normality for an empty lattice polytope $P$.

\begin{prop}
Let $P\subset M_{\mathbb{R}}$ be a lattice polytope of dimension $n$.
If 
$$
\mbox{\rm int}(nP)\cap M = \emptyset,
$$
then $P$ is isomorphic to a basic $n$-simplex.
\end{prop}
This fact is well known.  We see that
$$
\mbox{$P$ is a basic $n$-simplex} \Longleftrightarrow e(P)=n.
$$

\bigskip
We state a weaker condition than $k$-normality.
The following two lemmas are given in \cite{OZ}.

\begin{lem}\label{l:3.1}
Let $P\subset M_{\mathbb{R}}$ be a lattice polytope of dimension $n$.
If there exists an integer $r$ with $1\le r\le n-1$ such that $\mbox{\rm int}(rP)\cap M=\emptyset$,
then for all $k\ge n-r$ we have
$$
(kP)\cap M +(P\cap M) =(k+1)P\cap M.
$$
\end{lem}

\begin{lem}\label{l:3.2}
Let $P\subset M_{\mathbb{R}}$ be an empty lattice $n$-simplex.  If
$$
\mbox{\rm int}(n-1)P\cap M=\emptyset,
$$
then $P$ is basic.  
\end{lem}

\begin{lem}\label{l:3.3}
Let $P\subset M_{\mathbb{R}}$ be an empty lattice $n$-simplex. 
Assume that all facets of $P$ are basic.  If 
$$
\mbox{\rm int}(rP)\cap M=\emptyset \quad \mbox{and $r\ge \frac{n}2$},
$$
then $P$ is basic.
\end{lem}
{\it Proof}.   Lemma~\ref{l:3.1} saiys that 
$$
(kP)\cap M +(P\cap M)= (k+1)P\cap M \quad \mbox{for $k\ge n-r$}.
$$
From $r\ge \frac{n}2$, we have $n-r\le \frac{n}2$, hence $(n-r)P\cap M = \partial((n-r)P)\cap M$
from the assumption $\mbox{int}(rP)\cap M=\emptyset$ and $r\ge \frac{n}2$.
Thus for an $m\in (n-r)P\cap M$ there is a facet $F\subset P$ such that $m\in (n-r)F$.
Since $F$ is basic from the assumption, $m$ is a sum of $n-r$ elements of $F\cap M$.
This means that for $k\le n-r$ we have
$$
(kP)\cap M=\overbrace{(P\cap M)+ \dots +(P\cap M)}^k=(kP)\cap M.
$$
Therefor, $P$ is normal.
A normal empty lattice $n$-simplex is basic. \hfill $\Box$.

\bigskip
We calculate empty depth of several empty lattice simplices.

\begin{expl}
Let $e_1, \dots, e_n$ be a basis of free abelian group $M$ of rank $n\ge3$.
Consider an empty lattice $n$-simplex
$$
\Delta_2^n:=\mbox{\rm conv}\{0, e_1,e_2, e_1+e_2+2e_3, e_4, \dots, e_n\}.
$$
This polytope is not normal.  Since
$$
e_1+\dots+e_n\in \mbox{\rm int}(n-1)\Delta_2^n\cap M \setminus (n-2)\Delta_2^n,
$$
we have $e(\Delta_2^n)=n-2$. Thus
$$
\sharp(\Delta_2^n \cap M) -v(\Delta_2^n)=n-1 =e(\Delta_2^n)+1.
$$
\end{expl}

\begin{expl}\label{ex:3.2}
For $n\ge4$ we consider another empty lattice $n$-simplex
$$
\Delta_{2,n}^n:=\mbox{\rm conv}\{0, e_1, \dots, e_{n-1}, e_1+\dots +e_{n-1}+2e_n\}.
$$
Since
$$
e_1+\dots +e_n\in \frac{n}2 \Delta_{2,n}^n,
$$
we have $e(\Delta_{2,n}^n)=[\frac{n}2]$.  If $n\ge5$, then
$$
\sharp(\Delta_{2,n}^n \cap M) -v(\Delta_{2,n}^n)=n-1 > [\frac{n}2]+1 =e(\Delta_{2,n}^n).
$$
\end{expl}

\begin{rmk}
Example~\ref{ex:3.2} implies that the inequaliyu {\rm (\ref{e:p1})} does not hold in general for $n\ge5$.
\end{rmk}

\section{$h^*$ polynomials}

For a lattice polytope $P\subset M_{\mathbb{R}}$ of dimension $n$ we consider a power series
$$
\varphi(P,t):= \sum_{k\ge0}\sharp ((kP)\cap M) t^k.
$$
It is known that it has a form
\begin{equation}\label{e:4.1}
\varphi(P,t) =\frac{h_0^*+h_1^*t+ \dots +h_n^*t^n}{(1-t)^{n+1}}
\end{equation}
as a rational function (\cite{St}).  Here $h_i^*\ge0$ and $h_0^*=1$.  Moreover, it satisfies
\begin{equation}\label{e:v1}
h_0^*+h_1^*+\dots +h_n^*=v(P).
\end{equation}

\begin{dfn}
The polynomial $\sum_i h_i^*t^i=(1-t)^{n+1}\varphi(P,t)$ is called {\it $h^*$ polynomial} of $P$,
and written as $h_P^*$. The degree of the polynomial is called {\it degree} of $P$, written as
$\deg P$.
\end{dfn}

\begin{rmk}
For the basic $n$-simplex $\Delta^n$, we have
$$
\varphi(\Delta^n,t)=\frac1{(1-t)^{n+1}}.
$$
Hence, $\deg(\Delta^n)=0$.
\end{rmk}

\begin{rmk}
For a lattice polytope $P$, we set
$$
\varphi^*(P,t):= \sum_{k\ge0} \sharp(\mbox{\rm int}(kP)\cap M) t^k.
$$
From Ehrhart's reciprocity Theorem, we have
\begin{equation}\label{e:4.2}
\varphi^*(P,t)=\frac{h_n^*t+\dots +h_0^*t^{n+1}}{(1-t)^{n+1}}.
\end{equation}
In particular, we have
$$
h_n^*=\sharp(\mbox{\rm int}P\cap M).
$$
From this we see
\begin{equation}\label{e:4.3}
\deg P=n-e(P).
\end{equation}
\end{rmk}

From this terminology Batyrev and Nill give results \cite{BN}.

\begin{thm}[Batyrev-Nill]\label{th:BN1}
For a lattice polytope $P\subset M_{\mathbb{R}}$ of dimension $n$, the following two are equivalent;
\begin{itemize}
\item[(1)] $\deg P\le1$,
\item[(2)] The equality 
$$
\sharp(P\cap M) =v(P) +n
$$
holds.
\end{itemize}
\end{thm}

\begin{thm}[Batyrev-Nill]\label{th:BN2}
Let $e_1, \dots, e_n$ be a basis of a free abelian group $M$ of rank $n$.
A lattice polytope $P\subset M_{\mathbb{R}}$ of dimension $n$ satisfies $deg P=1$ if and only if
\begin{eqnarray*}
P&=& \mbox{\rm conv}\{0, e_1, \dots, e_{n-1}, e_1+a_1e_n< \dots, e_{n-1}+a_{n-1}e_n\}, a_ne_n\}\quad(n\ge2)\\
& & \mbox{where $a_1\ge \dots \ge a_{n-1}\ge0$, $a_n\ge1$ and $a_1+\dots +a_n\ge2$}
\end{eqnarray*}
or
$$
P=\mbox{\rm conv}\{0, 2e_1,2e_2, e_3, \dots, e_n\} \quad (n\ge3).
$$
The first polytope is called {Laurence prism} and the second {exceptional}.
\end{thm}

\section{The inequality}

In order to prove the inequality (\ref{e:p1}) we need some estimates of terms of $h^*$ polynomial.
 For a lattice polytope $P\subset M_{\mathbb{R}}$ of dimension $n$, we have
 \begin{equation}\label{e:5.1}
 \sharp(P\cap M) = h_1^* +n+1
 \end{equation}
by comparing the terms of degree one in the equality (\ref{e:4.1}).
From the equality (\ref{e:4.2}) we have
\begin{equation}\label{e:5.2}
\sharp(\mbox{int}P\cap M) =h_n^*, \quad \sharp(\mbox{int}(2P)\cap M) = h_{n-1}^* +(n+1)h_n^*.
\end{equation}

Since $\sharp(P\cap M) =\sharp(\partial P)\cap M +\sharp(\mbox{int}\ P\cap M)$ and $\sharp(\partial P\cap M)\ge n+1$,
from the equalities (\ref{e:5.1}) and (\ref{e:5.2}) we have
\begin{equation}\label{e:5.3}
h_1^* \ge h_n^*.
\end{equation}

\subsection{The case of $e\ge1$}

Let $P\subset M_{\mathbb{R}}$ be a lattice polytope of dimension $n$ satisfying $e(P)=n-1$.
Then $\deg h_P^*=1$ and the equality (\ref{e:v1}) is
$$
1+h_1^*=v(P).
$$
From the equality (\ref{e:5.1}) we have
$$
\sharp(P\cap M)-v(P)=n.
$$
This equality is just theorem~\ref{th:BN1}~(2).

\begin{lem}[The case of $e=n-2$ and $n\ge3$]\label{lem:5.1}
Set $n\ge3$. Let $P\subset M_{\mathbb{R}}$ be a lattice polytope of dimension $n$ with $e(P)=n-2$.
Then we have
$$
\sharp(P\cap M) -v(P) \le e(P)+1.
$$
The equality holds if and only if $\sharp((n-1)P\cap M)=1$.
\end{lem}
{\it Proof}. $h^*$ polynomial of $P$ is
$$
h_P^*=1+h_1^*t +h_2^*t^2.
$$
Since $e(P)=n-2$ we have $h_2^*=\sharp(\mbox{int}(n-1)P\cap M)\ge1$.
From (\ref{e:5.1}) and
$$
1+h_1^*+h_2^*=v(P),
$$
we have
$$
\sharp(P\cap M) -v(P) = n-h_2^* = (n-2)+1 -(h_2^*-1)\le e(P)+1.
$$
The equality holds if and only if $h_2^*=1$. \hfill $\Box$

\bigskip

\begin{rmk}
When the equality holds, $P$ defines a Gorenstein toric Del Pezzo variety. Batyrev and Juny classified 
Gorenstein toric Del Pezzo varieties {\rm \cite{BJ}}.
\end{rmk}

\bigskip

We need a general theory from Bruns and Herzog \cite{BH}.

\begin{prop}[Bruns-Herzog]\label{p:BH}
Let $P\subset M_{\mathbb{R}}$ be a lattice polytope of dimension $n$.  Set $e(P)=e$ and $h_P^*=h_0^*+
h_1^*t+\dots +h_{n-e}^*t^{n-e}$.  Then for $0\le l\le \frac{n-e}2$ we have
$$
\sum_{i=0}^l h_{n-e-i}^*\ge \sum_{i=0}^l h_i^*.
$$
\end{prop}

\begin{lem}[The case of $e=n-3$ and $n\ge4$]\label{lem:5.2}
Set $n\ge4$. Let $P\subset M_{\mathbb{R}}$ be a lattice polytope of dimension $n$ with $e(P)=n-3$.
Then  we have
$$
\sharp(P\cap M) -v(P) \le e(P)+1
$$
unless $n\ge5$ and $P$ is an empty lattice $n$-simplex isomorphic to
$$
\Delta_{2,5}^n:=\mbox{\rm conv}\{0,e_1, \dots, e_4, e_1+\dots + e_4+2e_5, e_6, \dots, e_n\}.
$$
\end{lem}
{\it Proof}. $h^*$ polynomial of $P$ is
$$
h_P^*=1+h_1^*t +h_2^*t^2+h_3^*t^3.
$$
Since $e(P)=n-3$ we have $h_3^*=\sharp(\mbox{int}(n-2)P\cap M)\ge1$.
From (\ref{e:5.1}) and
$$
1+h_1^*+h_2^*+h_3^*=v(P),
$$
we have
$$
\sharp(P\cap M) -v(P) = n-h_2^* -h_3^*= (n-3)+1 -(h_2^*+h_3^*-2).
$$
From Proposition~\ref{p:BH} we have
$$
h_2^*+h_3^*\ge 1+h_1^*.
$$
If $h_1^*\ge1$, then the inequality (\ref{e:p1}) holds.

$h_1^*=0$ implies that $P$ is an empty lattice simplex.
Set $n=4$.  Then $e(P)=1$ and $\sharp((2P)\cap M)=h_3^*\ge1$.
By comparing therms of degree two of (\ref{e:p1}) we have
$$
\sharp(2P\cap M) =\binom{n+2}2 +\binom{n+1}1 h_1^* +h_2^*.
$$
If $P$ is an empty lattice simplex, then $h_1^*=0$ and
$$
\sharp(\partial(2P)\cap M)\ge \binom{n+2}2.
$$
Hence $h_2^*\ge1$ and the inequality (\ref{e:p1}) holds when $n=4$.  

When $n\ge5$ the exception to the inequality is the case that $h_1^*=h_2^*=0$ and $h_3^*=1$, that is,
$P$ is an empty lattice $n$-simplex with $v(P)=2$.  Since $e(P)=n-3$, $P$ is isomorphic to $\Delta_{2,5}^n$. 
 \hfill $\Box$

\subsection{The case of $e=0$}

\begin{expl}\label{ex:5}
Let $M$ be a free abelian group of rank $n$ with a basis $e_1, \dots, e_n$.
Set
$$
\Delta_*^n:= \mbox{\rm conv}\{e_1, \dots, e_n.-e_1-\dots -e_n\}.
$$
The origin is only lattice point in the interior of $\Delta_*^n$ and $v(\Delta_*^n)=n+1$.
Its $h^*$ polynomial is
$$
h_{\Delta_*^n}^* = 1+ t+ \dots +t^n.
$$
In particular, $h_1^*=\dots =h_n^*=1$ and $\sharp(\Delta_*^n\cap M)-v(\Delta_*^n)=(n+2)-(n+1)=1$.
\end{expl}

Let $P\subset M_{\mathbb{R}}$ be a lattice polytope of dimension $n$ with $\mbox{int}P\cap M\not=\emptyset$.
Then $h_1^*\ge h_n^*\ge1$ from the inequality (\ref{e:5.3}). 
Moreover from Hibi's lower bound theorem \cite{Hb} we have $h_i^*\ge h_1^*$ for $2\le i\le n-1$.

\bigskip

\begin{prop}\label{p:5}
 Let $P\subset M_{\mathbb{R}}$ be a lattice polytope of dimension $n$ with $\mbox{\rm int}P\cap M\not=\emptyset$.
Then we have
$$
\Delta(P)\ge n-1.
$$
Moreover,when $n\ge3$ the equality holds if and only if $P\cong \Delta_*^n$.
\end{prop}
{\it Proof}. When $n=1$ $\Delta(P)=0$. Set $n\ge2$.
Set $h_P^*=\sum_{i=0}^n h_i^*t^i$.  We note $h_1^*, \dots , h_n^*\ge1$. 
From $v(P)=1+h_1^*+\dots +h_n^*$ and $\sharp(P\cap M) =h_1^*+n+1$, we have
\begin{equation}\label{e:5.10}
v(P)-\sharp(P\cap M)=(h_2^*-1)+\dots +(h_n^*-1)-1\ge -1.
\end{equation}
Hence we have
$$
\Delta(P) = n+v(P)-\sharp(P\cap M) \ge n-1.
$$

The equality holds if and only if $h_2^*=\dots=h_n^*=1$.
In particular, $h_n^*=\sharp(\mbox{int}P\cap M)=1$.
If $n\ge3$, then we know $h_2^*\ge h_1^*$.  Since $h_1^*\ge h_n^*=1$, we have $h_1^*=1$, 
 hence, $v(P)=n+1$ and $\sharp(P\cap M)=n+2$.
This means $P\cong \Delta_*^n$. \hfill $\Box$

\section{The case of dimension three}

In this section we give the main result of Tabei\cite{Tb}.

\begin{thm}[Tabei]\label{th:Tb}
Let $P\subset M_{\mathbb{R}}$ be a lattice polytope of dimension three.
Then we have
$$
\sharp(P\cap M) -v(P)\le e(P)+1.
$$
The lattice polytope $P$ satisfying the equality is one of the following:
\begin{itemize}
\item[(1)] When $e(P)=2$, $P$ is a Laurence prism or exceptional. 
\item[(2)] When $e(P)=1$, 
 \begin{itemize}
 \item[(2-1)] 16 cones over lattice polygons with one lattice point in their interiors.
 \item[(2-2)] 15 polytopes those facets have no lattice points in their interiors.
 \end{itemize}
 \item[(3)] When $e(P)=0$, 
 $$
 P\cong \mbox{\rm conv}\{(1,0,0),(0,1,0),(0,0,1),(-1,-1,-1)\}.
 $$
\end{itemize}
\end{thm}


\begin{thebibliography}{99}

\bibitem{BJ}
\textsc{V.~Batyrev and D.~Juny},
Classification of Gorenstein Toric Del Pezzo Varieties in arbitrary dimension,
Moscow Math. J. {\bf 10}, 285--316 (2010).

\bibitem{BN}
\textsc{V.~Batyrev and B.~Nill},
Multiples of Lattice Polytopes without interior Lattice points, 
Moscow Math. J. {\bf 7}, 195--207 (2007).

\bibitem{BH}
\textsc{W.~Bruns and J.~Herzog},
Cohen-Macaulay rings, Cambridge University Press, 1997.

\bibitem{Fj}
\textsc{T.~Fujita},
Classification theories of polarized varieties, London Mathematics Society Lecture Note Series {\bf 155}, Cambridge University Press (1990).

\bibitem{Ft}
\textsc{W.~Fulton},
Introduction to toric varieties, Ann. of Math. Studies No.131, Princeton Univ. Press, 1993.

\bibitem{Hb}
\textsc{T.~Hibi},
A lower bound theorem for Ehrhart polynomials of convex polytopes,
Adv. Math. {\bf 105} (1994), 162--165.

\bibitem{Ks}
\textsc{A. M.~Kasprzyk},
Canonical toric Fano threefolds, Canad. J. Math. {\bf 62} (2010), no. 6, 1293--1309.

\bibitem{KN}
\textsc{A. M.~Kasprzyk, M.~Kreuzer and B.~Nill},
On the combinatorial classification of toric log del pezzo surfaces,
LMS Journal of Computation and Mathematics {\bf 13} (2010), 33--46.

\bibitem{Ko}
\textsc{R.J.~Koelman},
The number of moduli of families of curves on toric surfaces,
Proefschrift, Katholieke Universitet te Nijmengen, 1991.

\bibitem{Od}
\textsc{T.~Oda},
Convex Bodies and Algebraic Geometry, Ergebnisse der Math. {\bf 15} Springer-Verlag,
Berlin, Heidelberg, New York, London, Paris, Tokyo, 1988.

\bibitem{OZ}
\textsc{S.~Ogata and H.-L.~Zhao},
A Characterization of Gorenstein toric Fano $n$-folds with index $n$
and Fujita's conjecture. Far East J. Math. Sci. {\bf 94} (2014), 65--88.

\bibitem{Sm}
\textsc{R.P.~Stanley},
Decompositions  of rational convex polytopes,
Ann. Discrete  Math. {\bf 6} (1980), 333--342.

\bibitem{Tb}
\textsc{R.~Tabei},
Toric varieties with maximal regularity, Tohoku University Master these, 2022.

\end{thebibliography}
\end{document}